
\input amstex

\magnification=1200
\loadmsam
\loadmsbm
\loadeufm
\loadeusm
\UseAMSsymbols

\hsize=6.0truein
\hoffset=0.15truein
\vsize=9truein
\voffset=-0.2truein

\def\leftitem#1{\item{\hbox to\parindent{\enspace#1\hfill}}}

\def\boxit#1#2{\hbox{\vrule
	\vtop{%
	\vbox{\hrule\kern#1%
	\hbox{\kern#1#2\kern#1}}%
	\kern#1\hrule}%
	\vrule}}

\def\leaderfill{\leaders\hbox to 1em{\hss.\hss}\hfill}

\parskip=\medskipamount
\document

\input epsf

\centerline{\bf Automorphisms of the Complex of Curves}

\centerline{ Feng Luo}

\centerline{\it Dept. of Math., Rutgers University, New Brunswick, NJ 08903}

\centerline{\it e-mail: fluo\@math.rutgers.edu}

\S1. {\bf Introduction}

	Given a compact orientable surface $\Sigma = \Sigma_{g,n}$
of genus g with n boundary components (possibly n = 0), let $\Cal S(\Sigma)$
 be the set of isotopy classes of essential unoriented non-boundary 
parallel simple loops in $\Sigma$. Two classes in $\Cal S(\Sigma)$
are called \it disjoint \rm
 if they are distinct and have disjoint representatives. 
In [Hav], Harvey introduced the complex of curves $\Cal C(\Sigma)$
for $\Sigma$ as follows. The vertices of $\Cal C(\Sigma)$ 
are elements in $\Cal S(\Sigma)$ and the simplexes of $\Cal C(\Sigma)$ are 
$<\alpha_1, ..., \alpha_k >$ where $\alpha_i$  is disjoint from 
$\alpha_j$ for $i \neq j$. This complex encodes the asymptotic geometry
 of the Teichm\"uller space in analogy with Tits buildings for symmetric spaces. The mapping class group acts on the curve complex preserving the simplicial structure. A natural question one would like to ask is whether every automorphism of the curve complex is induced by a homeomorphism of the surface.

In 1989 (see [Iv]), Ivanov sketched a proof of the result that if the
genus of the surface $g$ is at least 2, then
any automorphism of the curve complex $\Cal C(\Sigma_{g, n})$ 
is induced by a homeomorphism of the surface.

The aim of the paper is to settle the automorphism problem for 
the rest of the surfaces. Our proof does not distinguish the case 
genus $ g \geq 2$ from the case $g \leq 1$. We have,

{\bf Theorem. } \it (a) If the dimension $3g+n-4$ of the curve complex
$\Cal C(\Sigma_{g,n})$ at least  one and  $(g,n) \neq (1,2)$,  
then any automorphism of $\Cal C(\Sigma_{g,n})$ is 
induced by a self-homeomorphism of the surface.

	(b) Any automorphism  of $\Cal C(\Sigma_{1,2})$ preserving the
set of vertices represented 
by separating loops is induced by a self-homeomorphism of the surface.

	(c) There is an automorphism of $\Cal C(\Sigma_{1,2})$
which is not induced by any homeomorphisms.
\rm

We remark that Korkmaz [Ko] has independently proved part (a) of the theorem for genus $g \leq  1$ using different methods.

	We use the induction on the dimension of the curve complex 
$\Cal C(\Sigma)$ to prove the theorem. The strategy behind the proof fits 
extremely well with Grothendieck's philosophy (see [Gr]) that in the 
hierarchy of surfaces of negative Euler numbers under inclusion, the 
 ``generators" are the 1-holed torus and 4-holed sphere and the  ``relators" 
are the 2-holed torus and 5-holed sphere. Indeed, the most difficult and 
crucial cases in the proof are the 2-holed torus and 5-holed sphere 
whose curve complexes have dimension one. The proof for these 
two specific surfaces $\Sigma$ depends on an extremely simple fact that
 given two distinct elements in $\Cal S(\Sigma)$, there is at most one element in $\Cal S(\Sigma)$ 
which
is disjoint from both of them (lemma 4.3). Our proof makes extensive use 
of the work of several other authors ([Bi], [De], [Hal], [Iv], 
[Li] and [Vi]). In particular the work of Harer on the homotopy 
type of the curve complex is essential to our approach.

The theorem is an analogy to a result of Tits that all automorphisms of 
Tits buildings are induced by the automorphisms of the corresponding 
algebraic groups. 

Let Mod($\Sigma$) be the mapping class group Home($\Sigma$)/Iso of the 
surface of negative Euler number. There is a natural homomorphism 
$\pi$ 
: Mod($\Sigma$) $\to $  Aut($\Cal C(\Sigma)$) sending the isotopy 
class of a homeomorphism 
to the induced map on the curve complex. The theorem shows that $\pi$
is an epimorphism except $\Sigma= \Sigma_{1,2}$. The image 
$\pi $(Mod($\Sigma_{1,2})$) is a subgroup of index 5 in Aut($\Cal C(\Sigma_{0,5})$). By the work 
of Birman [Bi] and Viro [Vi], the kernel of $\pi$ is known to be trivial 
unless $(g,n) = (1,1),(0,4),(1,2),(2,0)$. If $(g,n)$ = $(1,1), (1,2)$ and
$(2,0)$, then the kernel is $\bold Z_2$ generated by a hyperelliptic
involution. If $(g,n) = (0,4)$, then  ker($\pi$) $ \cong$  $\bold Z_2 + \bold Z_2$
and is generated by two hyperelliptic involutions (see figure 6). 
Thus the theorem gives a new characterization of the mapping class group.

The complex of curves also arises in the study of 3-manifolds and mapping class groups. This complex was considered by Harer ([Har], [Har1]) 
from homological point of view (with applications to the homology of the mapping class group). In particular, Harer determined the homotopy
 type of the curve complex ([Hal], theorem 3.5). Ivanov ([Iv], [Iv1], 
[Iv2]) used the curve complex to determine the structure of the 
mapping class group. Masur and Minsky ([MM]) showed that the
 curve complex is $\delta$-hyperbolic in Gromov's sense. And Hempel
([He]) used the curve complex for studying 3-manifolds. 
See also  [Th].

The paper is organized as follows. In \S2, we establish basic properties of 
$\Cal C(\Sigma)$. In particular, using a result of Harer on the homotopy type of $\Cal C(\Sigma)$, 
it is shown that the curve complexes are pairwise non-isomorphic unless
$\Cal C(\Sigma_{1,1}) \cong \Cal C(\Sigma_{0,4})$, $\Cal C(\Sigma_{1,2})
\cong \Cal C(\Sigma_{0,5})$ and $\Cal C(\Sigma_{0,6}) \cong \Cal C(\Sigma_{2,0})$.
This is an analogy to Patterson's theorem for Teichm\"uller spaces. 
Part (c) of the theorem follows easily from $\Cal C(\Sigma_{1,2}) \cong
\Cal C(\Sigma_{0,5})$.
In  \S3, we introduce a multiplicative structure on $\Cal S(\Sigma)$. In particular, 
we define a $(\bold QP^1, PSL(2, \bold Z))$ modular structure on $\Cal S(\Sigma)$ 
(definition 3.4). The  $PSL(2, \bold Z)$ modular structure is fundamental to 
our approach to the automorphism problem. In  \S4, we show that any 
automorphism of $\Cal C(\Sigma)$ takes two curves intersecting at one point 
(resp. two points of different signs) to two curves intersecting at one point 
(resp. two points of different signs). Finally, in \S5, we prove the main
theorem 
by showing that any automorphism of $\Cal C(\Sigma)$ preserving the multiplicative 
structure is induced by a homeomorphism of the surface. To achieve th
is, we make extensive use of the modular structure (lemma 3.1).

\S2. {\bf Preliminaries}

2.1. {\it  Notations and conventions}

	We work in the piecewise linear category. All surfaces are 
oriented, connected and have negative Euler number. The isotopy class of a 1-dimensional submanifold $s$ is
 denoted by $[s]$. The group of homeomorphisms (resp. 
orientation preserving homeomorphisms) of $\Sigma$ is denoted by Home($\Sigma$) 
(resp. Home$^+(\Sigma)$). The group Home($\Sigma$) acts on $\Cal S(\Sigma)$ as follows:
$ h([a]) = [h(a)]$  where $h \in $ Home($\Sigma$) and $[a] \in $ $\Cal S(\Sigma)$. 
Given $\alpha, \beta \in $ $\Cal S(\Sigma)$, the
\it geometric intersection number \rm $I(\alpha, \beta)$ between the 
two classes is defined to be min\{$|a \cap b| | a \in \alpha, b \in \beta$\}.
If $F$ is a function defined on $\Cal S(\Sigma)$, we shall use $F(a)$ to denote 
$F([a])$ where  $[a] \in $ $\Cal S(\Sigma)$. In particular, if $a \in \alpha$, $b \in 
\beta$, then  $I(a, b) = I(a, \beta) =I(\alpha, b) = I(\alpha, \beta)$.
We shall use $\alpha \cap \beta = \emptyset$ to denote two disjoint elements 
 $\alpha$ and $\beta$, i.e., $I(\alpha, \beta) =0$ and $\alpha \neq \beta$.
If two elements $\alpha$ and $\beta$ satisfies $I(\alpha, \beta) \neq 0$, we
say that they \it intersect \rm and denote them by $\alpha \cap \beta
\neq \emptyset$. We use $\alpha \perp \beta$ to denote the relation 
$I(\alpha , \beta) = 1$.  And we use  $\alpha \perp_0 \beta$ to
denote two elements $\alpha$ and $\beta$ so that $I(\alpha, \beta) = 2$
 and their algebraic intersection number is zero.

	A subsurface  $\Sigma'$ in $\Sigma$  is called  \it
incompressible  \rm if the inclusion map $i: \Sigma' \to \Sigma$
induces a monomorphism in the fundamental groups. It is well known that 
 $\Sigma'$ is incompressible if and only if each component of 
$\partial \Sigma'$  is essential in  $\Sigma$. Assume that $\Sigma'$
is incompressible in $\Sigma$. Then the map $i_*: \Cal S(\Sigma')
\to $ $\Cal S(\Sigma)$ sending $[a]$ to $[i(a)]$
is injective so that $\alpha \cap \beta =\emptyset$, $\alpha \perp
\beta$, or $\alpha \perp_0 \beta$ in $\Cal S(\Sigma')$ if
and only if their images under $i_*$ satisfy the same relations.
 Due to this property, we shall identify  $\Cal S(\Sigma')$  with the 
subset  $i_*(\Cal S(\Sigma'))$.  An element $\alpha \in $ $\Cal S(\Sigma)$
is said to be in $\Sigma'$ if  $\alpha \in i_*(\Cal S(\Sigma'))$.
We say an element $\alpha$  in $\Cal S(\Sigma)$ \it
 decomposes \rm  $\Sigma$ into two subsurfaces 
 $\Sigma'$ and $\Sigma''$ if $\Sigma = \Sigma' \cup \Sigma''$ and
$\Sigma' \cap \Sigma'' \in \alpha$. If a class $\alpha \in \Cal S(\Sigma)$
decomposes the surface into a $\Sigma_{0,3}$ and $\Sigma'$, we say
$\alpha$ is a \it boundary class. \rm
A class $\alpha \in $ $\Cal S(\Sigma)$ is called \it separating \rm
 if it has a representative which is a separating loop.

Given a submanifold $s$, we use $ N(s)$ to denote a small regular 
neighborhood of $s$. We use $int(X)$ to denote the interior of a surface
$X$. The symbol $\cong$ is used to denote the homeomorphisms between surfaces, the isomorphisms between simplicial complexes, and isotopy.

Simple loops on surfaces will be denoted by small letters $a, b,$...,$ x, y, z$
and isotopy classes will be denoted by Greek letters $\alpha, \beta$, $\gamma$
etc.

2.2. {\it Basic properties of the curve complex }

	The homotopy type of the curve complex $\Cal C(\Sigma)$ was determined by 
Harer ([Hal], theorem 3.5).

{\bf Theorem (Harer)}.  \it
The curve complex  $\Cal C(\Sigma_{g,n})$ is homotopic to a wedge of
spheres of dimension $r$ where (i) $r=2g+n -3$, if $g>0$ and $n>0$,
(ii) $r = 2g-2$, if $n=0$ and (iii) $r= n-4$, if $g=0$. \rm

	A simplex in $\Cal C(\Sigma_{g,n})$ of maximal dimension
$3g+n-4$ is called a \it Fenchel-Nielsen system \rm
 (or a pant-decomposition). The following lemma is an easy 
consequence of Harer's theorem and Birman and Viro's work on the 
hyperelliptic involutions. The lemma is an analogous to a 
result of Patterson [Pa] 
that the Teichm\"uller spaces are pairwise nonisomorphic except
$T_{1,1} \cong T_{0,4}$, $T_{1,2} \cong T_{0,5}$, and
$T_{2,0} \cong T_{0,6}$.
Note that since $\Cal C(\Sigma_{1,1})$ and $\Cal C(\Sigma_{0,4})$
are zero-dimensional, by an isomorphism between them we mean a bijection 
 $\phi$ from $\Cal C(\Sigma_{1,1})$  to $\Cal C(\Sigma_{0,4})$
respecting the relations $\perp$ and $\perp_0$,
i.e., $\alpha \perp \beta$ if and only if 
$\phi(\alpha) \perp_0 \phi(\beta)$.

{\bf Lemma 2.1.}
\it
(a) $\Cal C(\Sigma_{2,0}) \cong \Cal C(\Sigma_{0,6})$,
$\Cal C(\Sigma_{1,2}) \cong \Cal C(\Sigma_{0,5})$ and
$(\Cal C(\Sigma_{1,1}), \perp) \cong (\Cal C(\Sigma_{0,4}),\perp_0)$.

(b) If $(g,n) \neq (g',n')$ and  $\Cal C(\Sigma_{g,n})$
is  not one of the six complexes above, then the curve 
complexes $\Cal C(\Sigma_{g,n})$ and $\Cal C(\Sigma_{g', n'})$ are
not isomorphic.
\rm

{\it Proof.} To show (a), let us first consider $\Cal C(\Sigma_{0,6})
\cong$  $\Cal C(\Sigma_{2,0})$.
We construct a bijection from $\Cal S(\Sigma_{0,6})$ to $\Cal S(\Sigma_{2,0})$
preserving the disjointness 
as follows. Let $r: \Sigma_{2,0} \to \Sigma_{2,0}$ be a hyperelliptic 
involution. Then  $r(\alpha) = \alpha$ for all $\alpha \in 
\Cal S(\Sigma_{2,0})$ by a result of Birman ([Bi]) and Viro ([Vi]). 
Indeed, $r$ commutes with 
all Dehn twists up to isotopy. Let $P: \Sigma_{2,0} \to \Sigma_{2,0}/ r
\cong S^2$ be the quotient map which is a branched covering branched 
over a six-point set $B$. Define $P^*: \Cal S(S^2 - int(N(B)))
\to \Cal S(\Sigma_{2,0})$ by sending $[a]$ to $[b]$ where $b$ is a 
component of  $P^{-1}(a)$. Then  $P^*$ is a bijection preserving disjointness. 
Now, $S^2 - int(N(B)) \cong \Sigma_{0,6}$. Thus $\Cal C(\Sigma_{0,6})
\cong \Cal C(\Sigma_{2,0})$.

For $\Sigma_{1,2}$, take a non-separating $r$-invariant simple loop $s$ 
and let $\Sigma_{1,2} = \Sigma_{2,0} - int(N(s))$. Then $P(\Sigma_{1,2})$
 is a disc with 4-cone points $B_4$ of order 2. 
Let $\Sigma_{0,5}$ be $P(\Sigma_{1,2}) - int(N(B_4))$. 
Then $P^*|_{\Cal S(\Sigma_{0,5})}$ is a bijection from 
$\Cal S(\Sigma_{0,5})$ onto  $\Cal S(\Sigma_{1,2})$ preserving disjointness.
 Finally, identify $\Sigma_{1,1}$ an $r$-invariant subsurface of 
$\Sigma_{2,0}$. Then $P(\Sigma_{1,1})$  is a disc with three 
cone points of order 2. The same argument shows that the restriction of 
$P^*$ gives a bijection between $\Cal S(\Sigma_{0,4})$ and 
$\Cal S(\Sigma_{1,1})$ 
which respects the relations $\perp_0$ and $\perp$.

To see (b), take $(g, n) \neq (g', n')$. Using Harer's theorem and counting 
the dimension of the curve complex, we conclude that that $\Cal 
C(\Sigma_{g,n})$
and  $\Cal C(\Sigma_{g', n'})$ are not isomorphic except possibly 
the following cases: (i) $(g',n') = (0,n')$ with $n' \geq 7$ and $(g,n)$ with 
$g \geq 1$,  and (ii) $(g,n) = (g,3) $ and $(g',n') = (g + 1,0)$. 

In  case (i), suppose otherwise that $\phi : \Cal S(\Sigma_{g,n})
\to \Cal S(\Sigma_{0,n'})$ is a bijection preserving disjointness. Since
 $g \geq 1$, take a non-separating class $\alpha \in 
\Cal S(\Sigma_{g,n})$. Then 
 $\phi(\alpha)$ must be a boundary class, i.e., it decomposes $\Sigma_{0,n'}$
into an $\Sigma_{0,3}$ and $\Sigma'$. To see this, for any two classes 
 $\beta$ and $\gamma$ disjoint from $\alpha$, there exists a sequence
$\alpha_1 = \beta, \alpha_2, ..., \alpha_k = \gamma$ so that
$\alpha_i \cap \alpha = \emptyset$ and $\alpha_i \cap \alpha_{i+1}
\neq \emptyset$. 
However, if $\phi(\alpha)$ is not a boundary class, there exist two classes 
 $\beta '$ and $\gamma'$ disjoint from $\phi(\alpha)$  which cannot be
 joint by such a sequence. Since $g \geq 1$,  there exists a maximal 
dimension simplex $<\alpha_1, ..., \alpha_k>$ in $\Cal C(\Sigma_{g,n})$
so that each $\alpha_i$ is non-separating. Its image under $\phi$ 
is a maximal dimension simplex in $\Cal C(\Sigma_{0,n'})$ so that each
vertex is a boundary class. This is impossible unless
$n'= 4,5, $ or $6$.

In case (ii), suppose otherwise that $\phi: \Cal C(\Sigma_{g,3})
\to \Cal C(\Sigma_{g+1, 0})$ is an isomorphism where $g \geq 1$. Take 
a nonseparating class $\alpha \in \Cal S(\Sigma_{g, 3})$  and consider 
its image under $\phi$.
 Since there are no boundary classes in $\Cal S(\Sigma_{g+1, 0})$, the
image $\phi(\alpha)$  must be non-separating by the same argument as before. 
By considering the classes disjoint from $\alpha$, we obtain the following 
isomorphism  $\Cal C(\Sigma_{g-1, 5})  \cong \Cal C(\Sigma_{g,2})$.
By the result just proved above, this shows $g=1$, i.e., we have
$\Cal C(\Sigma_{1,3}) \cong \Cal C(\Sigma_{2,0})$. 
But by part (a), we have $\Cal C(\Sigma_{2,0}) \cong \Cal C(\Sigma_{0,6})$.
Thus we obtain $\Cal C(\Sigma_{1,3}) \cong \Cal C(\Sigma_{0,6})$.
This contradicts the conclusion of  case (i).
$\square$

\it Remark. \rm  The maximal dimension of those simplexes
$<\alpha_1, ..., \alpha_k>$ in $\Cal C(\Sigma_{g,n})$ so that each 
$\alpha_i$ is separating is $2g+n-4$. 

{\it Proof of part (c) of the main theorem}. 
Since $\Cal C(\Sigma_{1,2})$ is isomorphic to $\Cal C(\Sigma_{0,5})$
and Home($\Sigma_{0,5})$ acts transitively on $\Cal S(\Sigma_{0,5})$, 
the automorphism group of $\Cal C(\Sigma_{1,2})$ acts transitively on
$\Cal S(\Sigma_{1,2})$.  In particular, there is an automorphism of
$\Cal C(\Sigma_{1,2})$ which sends a separating class to a
non-separating class. $\square$

{\bf Lemma 2.2.} \it Suppose $3g+n \geq 5$ and $(g,n) \neq (1,2)$.
If  $\phi:  \Cal S(\Sigma_{g,n}) \to \Cal S(\Sigma_{g,n})$
is a bijection preserving disjointness,
then  $\phi$  preserves the separating classes. \rm

{\it Proof.}  Suppose otherwise that $\alpha$  is non-separating and 
$\phi(\alpha)$  is separating. Then by the same argument as in the proof of 
lemma 2.1, $\phi(\alpha)$  is a boundary class. 
By considering the isotopy classes disjoint from $\alpha$ and $\phi(\alpha)$
respectively,  we obtain the following isomorphisms
$\Cal C(\Sigma_{g-1, n+2}) \cong \Cal C(\Sigma_{g, n-1})$. By lemma 2.1, 
we conclude that $(g,n) = (1,2)$ or $(1,3)$.  By the assumption, thus  $(g,n) = (1,3)$.  Extend $\alpha$  to a Fenchel-Nielsen system $\{\alpha, \beta, \gamma\}$
 so that both $\beta$ and $\gamma$ are non-separating. Thus 
$\{ \phi(\alpha), \phi(\beta), \phi(\gamma)\}$ is again a Fenchel-Nielsen 
system in $\Sigma_{1,3}$  where $\phi(\alpha)$ is a boundary class.
Say, $\phi(\alpha)$   bounds a subsurface $\Sigma_{1,2}$.  Since any 
Fenchel-Nielsen system on $\Sigma_{1,2}$  contains a non-separating element, 
thus one of the element $\phi(\beta)$ or $\phi(\gamma)$ is non-separating in the 
subsurface $\Sigma_{1,2}$. Say  $\phi(\beta)$ is non-separating. Find a class 
 $\delta$ in $\Sigma_{1,2}$ which is disjoint from  $\phi(\beta)$
so that $\delta$  bounds a $\Sigma_{1,1}$ in $\Sigma_{1,2}$. Thus 
 $\delta$ decomposes $\Sigma_{1,3}$ into a union of $\Sigma_{1,1}$ and $\Sigma_{0,4}$. By lemma 2.1, $\phi^{-1}(\delta)$ bounds a $\Sigma_{1,1}$ in $\Sigma_{1,3}$.
Thus we have a Fenchel-Nielsen system $\{\alpha, \beta, \phi^{-1}(\delta)\}$
so that $\alpha$ and $\beta$ are non-separating and $\phi^{-1}(\delta)$
 bounds $\Sigma_{1,1}$.  This is absurd. Thus the lemma follows. $\square$

The following generalizes an observation of Ivanov which he proved for genus at least 2.

{\bf Lemma 2.3.} \it If $\phi: \Cal S(\Sigma) \to \Cal S(\Sigma)$ is a 
bijection preserving disjointness
so that in the case $\Sigma = \Sigma_{1,2}$, $\phi$ preserves the separating
classes, then for any $\alpha \in \Cal S(\Sigma)$ there exists $h \in $ Home($\Sigma)$ so that $h(\alpha) = \phi(\alpha)$. \rm

{\it Proof}.  If $\Sigma = \Sigma_{g,n}$ satisfies $3g+n \leq 5$, 
then the lemma is evident. Assume now that $3g+n \geq 6$. By lemma 2.2, 
we may assume further that
$\alpha$ is separating and is not a boundary class. 
Take $a  \in \alpha$ and $b \in  \phi(\alpha)$.
By lemma 2.2, we have $\Sigma -int(N(a)) = \Sigma_{g_1, n_1} \cup
\Sigma_{g_2, n_2}$ and $\Sigma - int(N(b)) = \Sigma_{g'_1, n'_1} \cup
\Sigma_{g'_2, n'_2}$ so that the curve complexes
$\Cal C(\Sigma_{g_i, n_i})$ and $\Cal C(\Sigma_{g'_i, n'_i})$ are both
non-empty. The goal is to show that the two decompositions are homeomorphic.
To this end, we note first that $g_1 + g_2 = g = g_1' + g_2'$ and
$n_1 + n_2 =n = n_1' + n_2' $. Second, the bijection $\phi$ sends
the pair $\{ \Cal S(\Sigma_{g_1, n_1}), \Cal S(\Sigma_{g_2, n_2})\}$ to
$\{ \Cal S(\Sigma_{g_1', n_1'}), \Cal S(\Sigma_{g_2', n_2'})$\} by the
same argument as in the proof of lemma 2.1. Assume without loss of
generality that $\phi(\Cal C(\Sigma_{g_i, n_i})) = \Cal C(\Sigma_{g_i', n_i'})$.
It remains to show that $(g_i, n_i) = (g_i', n_i')$ for $i=1,2$ in order
to finish the proof.

	By lemma 2.1, we obtain $(g_i, n_i) = (g_i', n_i')$ except the 
following three decompositions which need to be checked specifically. Namely (i)
$(g_1, n_1) =(g_2', n_2')=(1,1)$ and $(g_2, n_2) = (g_1', n_1') =(0,4)$;
(ii) $(g_1, n_1) = (g_2', n_2') =(0,5)$ and $(g_2, n_2) = (g'_1, n'_1)
=(1,2)$; and (iii) $(g_1, n_1) = (1,1)$, $(g_1', n_1') = (0,4)$,
$ (g_2, n_2) = (0,5)$ and $(g_2', n_2') = (1,2)$.
None of these three cases occurs due to lemma 2.2. Indeed, if there 
were $\alpha \in \Cal S(\Sigma)$ decomposing the surface $\Sigma$ into 
a genus 0 and a genus 1 subsurfaces and $\phi$ interchanges the two 
subsurfaces, then $\phi$ would send a non-separating class to a 
separating class.
$\square$.

Given an incompressible subsurface $\Sigma'$ in a surface $\Sigma$, 
then each non-separating
 class in $\Cal S(\Sigma')$ is again non-separating in $\Cal S(\Sigma)$. But separating 
classes in $\Cal S(\Sigma')$ may become nonseparating in $\Cal S(\Sigma)$. However, if 
 $\Sigma_{1,2}$ is an incompressible subsurface in a surface $\Sigma$, then 
each separating class in $\Sigma_{1,2}$ remains separating in $\Sigma$. 
Thus given  an incompressible subsurface $\Sigma_{1,2}$, the inclusion map from 
$\Sigma_{1,2}$ to $\Sigma$
preserves both separating and nonseparating classes. Combining with lemma 2.3, 
we have,

{\bf Corollary 2.4.} \it Suppose the dimension of $\Cal C(\Sigma)$ is
at least 2 and $\Sigma_{1,2}$ is an incompressible subsurface in $\Sigma$.
 If  $\phi: \Cal S(\Sigma) \to \Cal S(\Sigma)$ is  a bijection preserving the 
disjointness so that $\phi(\Cal S(\Sigma_{1,2})) 
= \Cal S( \Sigma_{1,2})$, then
$\phi|_{\Cal S(\Sigma_{1,2})}$ preserves the separating classes.
\rm

{\bf  Corollary 2.5}. \it Suppose $\alpha_1 \in \Cal S(\Sigma)$ decomposes
$\Sigma$ into a union of $\Sigma'$ and $\Sigma''$ so that $\Sigma'
\cong \Sigma_{1,1}$ or $\Sigma_{0,4}$ and suppose $\alpha_2 \in
\Cal S(\Sigma')$. Then given any bijection  $\phi$
of $\Cal S(\Sigma)$ preserving the
disjointness and the separating classes, there is a homeomorphism
$h$ of the surface $\Sigma$ so that $h(\alpha_i) = \phi(\alpha_i)$ for
$i=1,2$. \rm

\it Proof. \rm  First, we claim that there is a homeomorphism $h_1$
 of the surface $\Sigma$ so that $h_1(\alpha_1) = \phi(\alpha_1)$ 
and $h_1(\Cal S(\Sigma')) = $ $\phi(\Cal S(\Sigma'))$.  Indeed, by 
lemma 2.3, we find $h_2 \in$
Home($\Sigma$) so that $h_2(\alpha_1) = \phi(\alpha_1)$. By the proof of
lemma 2.3, $h_2(\Cal S(\Sigma')) = \phi(\Cal S(\Sigma'))$ unless $\Sigma'' \cong
\Sigma_{1,1}$ or $\Sigma_{0,4}$. If  $h_2(\Cal S(\Sigma'))
=\Cal S(\Sigma'')$,  then $\Sigma' \cong \Sigma''$ because $\phi$ 
preserves separating classes. Now let $h_3$
be an involution of $\Sigma$ interchanging $\Sigma'$ and $\Sigma''$ and fixing
$\alpha_1$. Then $h_1 = h_3 \circ h_2$ is a required homeomorphism.

Second, let $a_1$ be the component of $\partial \Sigma'$ corresponding
to $\alpha_1$. The group of homeomorphisms of $\Sigma'$ leaving $a_1$
pointwise fixed acts transitively on $\Cal S(\Sigma')$. Thus, we
 find a homeomorphism $h_4$ of $\Sigma$ which is the identity map on
$\Sigma''$ so that $h_4(\alpha_2) = \phi(\alpha_2)$. 
The required homeomorphism $h = h_4 \circ h_1$.
$\square$

\S3. {\bf A ($\bold QP^1, PSL(2, \bold Z))$ Structure on $\Cal S(\Sigma)$}

3.1. Recall that surfaces are oriented. Suppose $x$ and $y$ are two
open arcs intersecting transversely at a point $p$ in a surface $\Sigma$.
Then the  \it resolution of the intersection point $p$ from $x$ to $y$ \rm
is defined as follows. Fix an orientation on $x$. Use the orientations on
the surface $\Sigma$ and $x$  to determine an orientation on $y$. Finally
 resolve the intersection according to the orientations (see figure 1(a)).
 This resolution is independent of the choice of orientations on $x$.
 Suppose now that $\alpha$ and $\beta$ are two elements in $\Cal S(\Sigma)$ with
$\alpha \perp \beta$ or $\alpha \perp_0 \beta$,  we define the
 multiplication $\alpha \beta$ as follows. Take $a \in \alpha$ and $b \in
\beta$  so that $|a \cap b| = I(\alpha, \beta)$.  Then $\alpha \beta$
is the isotopy class $[ab]$  where $ab$ is the simple loop obtained by
 resolving all intersection points in $a\cap b$ from a to b. See figure 1(a).
One checks easily that if $\alpha \perp \beta$ (resp. $\alpha \perp_0 \beta$)
then $\alpha \beta \in \Cal S(\Sigma)$ and $\alpha \beta \perp \alpha$, $\beta$
(resp. $\alpha \beta \perp_0 \alpha, \beta$).

3.2. Let $\hat \bold Q = \bold Q \cup \{\infty\}$. Two rational numbers
$p/q$ and $p'/q'$ satisfying $pq' - p'q = \pm 1$ are denoted by
$p/q \perp p'/q'$. The relation $(\hat \bold Q, \perp)$ is called the
\it modular configuration\rm.  A standard way of presenting the 
configuration is to consider $\hat \bold Q$ as a subset of the boundary
of the upper half plane $\bold H$ and to draw a hyperbolic geodesic ending
at $p/q$ and $p'/q'$ if $p/q \perp p'/q'$.  Figure 1(b) is the configuration
after a M\"obius transformation. It was known to Max Dehn [De] that both
$(\Cal S(\Sigma_{1,1}), \perp)$ and $(\Cal S(\Sigma_{0,4}), \perp_0)$
are isomorphic to the modular configuration, i.e., there exists a
bijection  $\pi$ between $\Cal S(\Sigma_{1,1})$ (resp. $\Cal S(\Sigma_{0,4}))$
and $\hat \bold Q$ so that $\alpha \perp \beta$ (resp. $\alpha \perp_0 \beta$)
 if and only if $\phi(\alpha) \perp \phi(\beta)$.  Furthermore, if
$\pi(\alpha) = p/q$ and $\pi(\beta) = p'/q'$ so that $p/q \perp p'/q'$,
then $\pi(\alpha \beta) = (p+\lambda q)/(p' + \lambda  q')$ where
$\lambda = pq'-p'q$.
Note that $(\pi(\alpha), \pi(\alpha \beta), \pi(\beta))$ determines
the right-hand orientation on the circle.

\midspace{0.1cm}
\centerline{\epsfbox{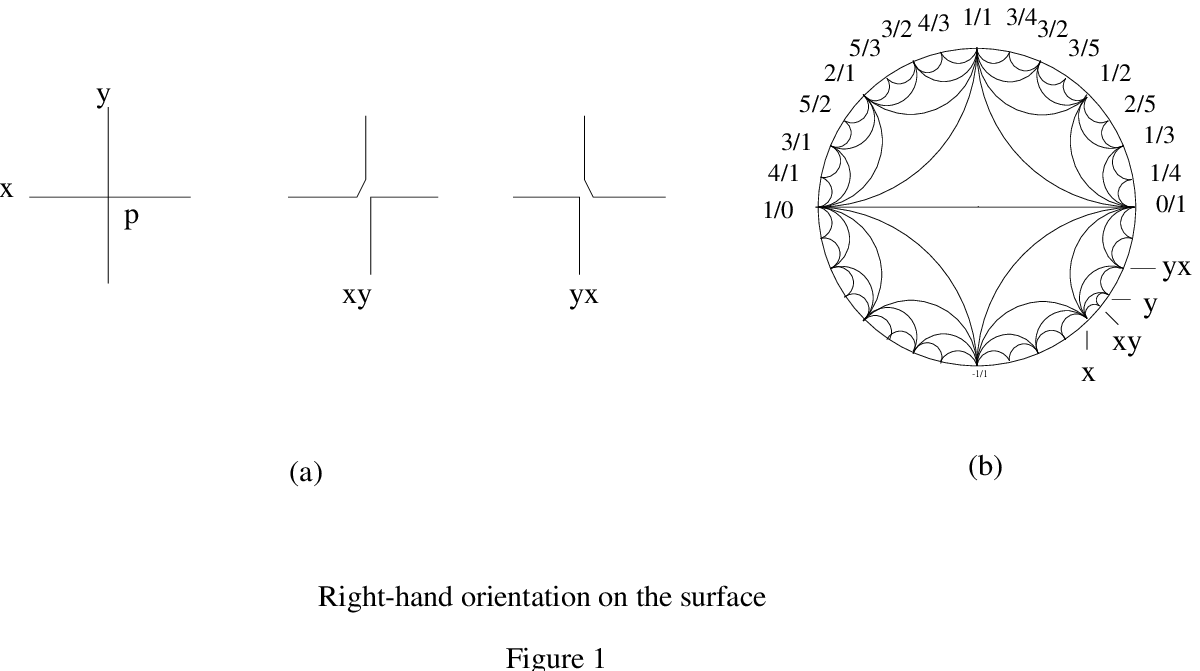}}
\midspace{0.1cm}

The following lemma is an easy consequence of the modular configuration.

{\bf Lemma 3.1.} \it (a) If $\phi : \Cal S(\Sigma_{1,1}) \to \Cal 
S(\Sigma_{1,1})$ (resp. $\Cal S(\Sigma_{0,4}) \to
\Cal S(\Sigma_{0,4}))$ is a bijection preserving the relation $\perp$ 
(resp. $\perp_0$), then $\phi$ is induced by a homeomorphism of the surface.

(b) Two elements $\alpha_1, \alpha_2 \in \hat \bold Q$
satisfy  $\alpha_1 \perp \alpha_2$  if and only if there are
two distinct elements $\gamma_1, \gamma_2$ so that $\gamma_i \perp \alpha_j$
and $\gamma_1$ and $\gamma_2$ are not related by $\perp$. Furthermore,
$\{ \gamma_1, \gamma_2 \} =\{ \alpha_1 \alpha_2, \alpha_2 \alpha_1\}$.
\rm

3.3. We begin by introducing some notations. If $\alpha \perp \beta$
or $\alpha \perp_0 \beta$,  we denote it by 
$\alpha \top \beta$.  Given a subset $\Cal X \subset$ $\Cal S(\Sigma)$, let $\Cal X_{\infty}
= \cup_{n=0}^{\infty} \Cal X_n$ where $\Cal X_{0} = \Cal X$,
and $\Cal X_{n+1} = \Cal X_n \cup \{\alpha | \alpha = \beta \gamma$, where
$\beta \top \gamma,$ and $\beta, \gamma$, $\gamma \beta$ are in
$\Cal X_n$\}. If $\Cal X_{\infty} = $ $\Cal S(\Sigma)$, we say that $\Cal X$ \it
generates \rm $\Cal S(\Sigma)$. For instance, the three-element set
$\{\alpha, \beta, \alpha \beta\}$  
generates the sets $\Cal S(\Sigma_{1,1})$ and $\Cal S(\Sigma_{0,4})$.
The following lemma is motivated by the proof of lemma 2 in [Li]. See 
also [Lu].

{\bf Lemma 3.2.} \it  Suppose  $\{\alpha_1, ... \alpha_k\}$ are pairwise
disjoint elements in $\Cal S(\Sigma)$ and $\alpha \in $ $\Cal S(\Sigma)$ so that $I(\alpha, \alpha_1)
\geq 2$ and $\alpha$ and $\alpha_1$ are not related by $\perp_0$.
Then $\alpha = \beta_1 \beta_2$ where $\beta_1 \top \beta_2$ so that
$I(\beta_i, \alpha_1) < I(\alpha, \alpha_1), I(\beta_2 \beta_1, \alpha_1)
< I(\alpha, \alpha_1)$, $I(\beta_i, \alpha_j) \leq I(\alpha, \alpha_j)$,
and $I(\beta_2 \beta_1, \alpha_j) \leq I(\alpha, \alpha_j)$ for $i=1,2$
and $j \geq 2$. In particularm $\Cal S(\Sigma)$ is generated by the set
$\Cal G =\{ \alpha \in \Cal S(\Sigma) |$ for each $i$, either $\alpha
\top \alpha_i$ or $\alpha \cap \alpha_i = \emptyset$\}. \rm

{\it Proof.}  Take $a \in \alpha$ and $a_i \in \alpha_i$ so that
$|a \cap a_i| =I(\alpha, \alpha_i)$ and $a_i \cap a_j = \emptyset$. Now
consider the following two cases.

{\it Case 1.} There exist two points $p,q \in a \cap a_1$ which are 
adjacent in $a_1$ so that they have the same intersection sign. See  
figure 2.

\midspace{0.1cm}
\centerline{\epsfbox{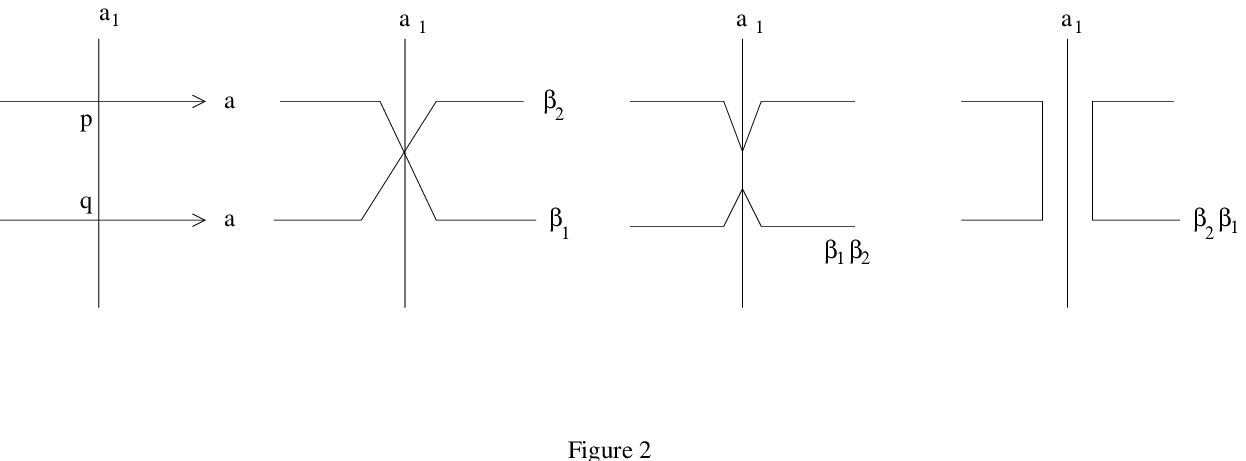}}
\midspace{0.1cm}

Assuming that the surface has the right-hand orientation, we take $\beta_1$
and $\beta_2$  as indicated. Then $\beta_1 \perp \beta_2$, and $\alpha
= \beta_1 \beta_2$. We verify the required conditions for $\beta_1$ and
$\beta_2$  as in figure 2. If the surface has the 
left-hand orientation, we interchange $\beta_1$ and $\beta_2$.

{\it Case 2.} If case 1 does not occur, 
then there are three points $p,q$ and $r$ in  $a \cap a_1$
which are adjacent in $a_1$ so that their intersection signs 
alternate. See figure 3.
Fix an orientation on  $a$. Assume without loss of generality that 
the arc in $a$ from $p$ to  $q$ does not contain the point $r$. 
If the surface has the right-hand orientation, we choose $\beta_1$ and
$\beta_2$ as in figure 3. Since $|a \cap a_1| = I(\alpha, \alpha_1)$,
$\beta_1$ and $\beta_2$ are both in $\Cal S(\Sigma)$ and $\beta_1 \perp_0 \beta_2$.
We have $\alpha = \beta_1 \beta_2$.
The required conditions for $\beta_i$'s are verified as in figure 3. 
If the surface has the left-hand orientation, we interchange $\beta_1$ and
$\beta_2$. $\square$

\midspace{0.1cm}
\centerline{\epsfbox{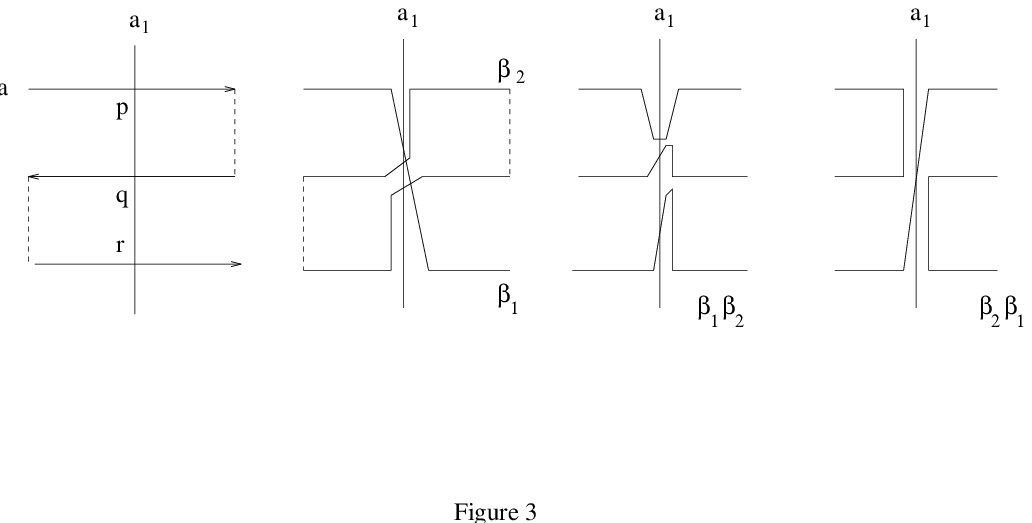}}
\midspace{0.1cm}

\it Remark. \rm   A stronger version of the lemma still holds. See [Lu]
lemma 7.

{\bf Corollary 3.3.} \it Under the same assumption as in lemma 3.2, suppose 
$\phi$ and $\psi$ are two bijections of $\Cal S(\Sigma)$ satisfying the 
following conditions:

(1) Both $\phi$ and $\psi$  preserve the disjointness and relations
$\perp$ and $\perp_0$.

(2) If $\alpha \top \beta$, then  $\{ \phi(\alpha \beta), \phi(\beta \alpha)\}
=\{ \phi(\alpha) \phi(\beta), \phi(\beta) \phi(\alpha)\}$ and
$\{ \psi(\alpha \beta),$$ \psi(\beta \alpha)\}$\newline
$ =\{ \psi(\alpha) \psi(\beta),$$ \psi(\beta) \psi(\alpha)\}$.

(3)  $\phi|_{\Cal G} = \phi|_{\Cal G}.$

Then $\phi = \psi$. \rm

3.4. We have mentioned in several places the notion of modular structure on a discrete set. Here  is a formal definition after Thurston's 
geometric structures on manifolds.

{\bf Definition.} A \it modular structure \rm on a discrete set X is a maximal 
collection of charts $\{ (U_i, \phi_i) | i \in I\}$ where $\phi_i: U_i \to
\bold QP^1$ is injective so that the following three conditions 
are satisfied:

(1) The union of the domains of the charts covers X, i.e., $\cup_{i \in I}
U_i = X$.

(2) The transition functions  $\phi_i \phi_j^{-1}$  are restrictions of 
elements in  $PSL(2, \bold Z).$

A modular structure on a set $X$ is called \it compact \rm if the
following condition holds,

(3) The automorphism group of the structure $(X, \{(U_i, \phi_i)\})$
acts on X with finite orbits.

The last condition seems to be crucial. Examples of modular structure are 
$\Cal S(\Sigma)$ and the set of all Fenchel-Nielsen systems (see [Lu2]).

{\bf Lemma 3.4.} \it  If $\Sigma$ is an oriented surface with $\Cal S(\Sigma)$ $\neq \emptyset$,
then $\Cal S(\Sigma)$ has a natural modular structure invariant under the 
action of the orientation preserving mapping class group. \rm

In fact, as a consequence of the main theorem of the paper, one sees that 
the automorphism group of the modular structure on $\Cal S(\Sigma)$ is the 
orientation preserving mapping class group for all surfaces.

{\it Proof. } If the dimension of $\Cal C(\Sigma)$ is zero, then the surfaces are 
$\Sigma_{1,1}$, $\Sigma_{0,4}$ or $\Sigma_{1,0}$.
The result follows by the proof of lemma 2.1. Fix a standard oriented 
1-holed torus $\Sigma_{1,1}$ and an identification between 
$\Cal S(\Sigma_{1,1})$ and $ \bold QP^1$. If the dimension of the complex 
$\Cal C(\Sigma)$ is at least one, then any element in $\Cal S(\Sigma)$ lies in an incompressible 
subsurface $\Sigma'$  homeomorphic to either $\Sigma_{1,1}$ or
$\Sigma_{0,4}$. Assume the subsurface has the induced orientation. 
Then the charts are  $(\Cal S(\Sigma'), \phi)$ where $\phi: \Cal S(\Sigma')
\to \Cal S(\Sigma_{1,1})$ is a bijection produced in the proof 
of lemma 2.1 so that $\phi$  respects the orientations. Extends 
these charts to be the maximal collection. One checks easily that all 
conditions are satisfied. $\square$ 

\S4. {\bf  A Basic Property of the Automorphisms of $\Cal S(\Sigma)$}

The aim of the section is to prove the following proposition.

{\bf Proposition.} \it  Suppose $3g +n \geq 5$ and $\phi$ : $\Cal S(\Sigma)$ $\to $ $\Cal S(\Sigma)$
is a bijection preserving disjointness and the separating classes. 
Then $\phi$  preserves the relations  $\perp$ and $\perp_0$ 
in $\Cal S(\Sigma)$. 
Furthermore, if $\alpha \top \beta$,  then 
 $\{ \phi(\alpha \beta), \phi(\beta \alpha)\}
=\{ \phi(\alpha) \phi(\beta), \phi(\beta) \phi(\alpha)\}$.\rm

{\it Proof.} We use induction on $|\Sigma_{g,n}| = 3g+n$.  The main step 
is in the case where $|\Sigma| =1$, i.e., $\Sigma = \Sigma_{0,5}$ (case 1)
 and $\Sigma_{1,2}$ (case 2).

\midspace{0.1cm}
\centerline{\epsfbox{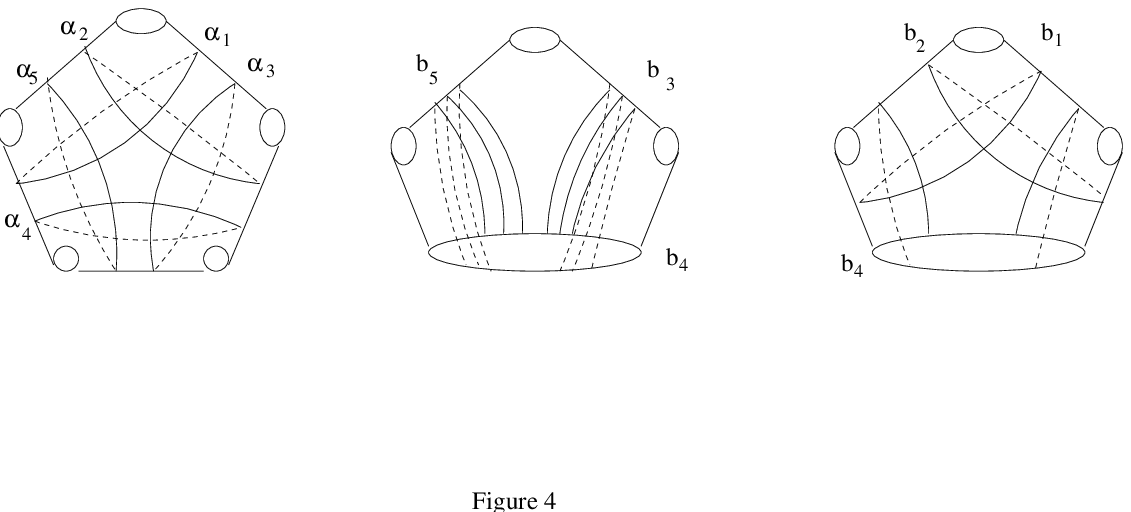}}
\midspace{0.1cm}

{\it Case 1}. The surface is  $\Sigma_{0,5}$.
We first show that $\phi$ preserves the relation $\perp_0$. 
To this end, take two isotopy classes $\alpha_1$ and $\alpha_2$ so that
$\alpha_1 \perp_0 \alpha_2$.  To show that 
$\phi(\alpha_1) \perp_0 \phi(\alpha_2)$, we extend $\{\alpha_1, \alpha_2\}$
to a ``pentagon" $\{\alpha_1, ..., \alpha_5\}$ where $\alpha_i \perp_0
\alpha_{i+1}$ and $\alpha_i \cap \alpha_{i+2} = \emptyset$ (indices
$i$ are counted mod 5) as in figure 4. Here we have used a simple
fact that any two pairs of isotopy classes $(\alpha, \beta)$ with
$\alpha \perp_0 \beta$ in $\Cal S(\Sigma_{0,5})$ are related
by a homeomorphism of the surface. Indeed,  if we take $a \in
\alpha$ and $b \in \beta$ with $|a \cap b| = 2$, then the
regular neighborhood $N(a \cup b)$ is an incompressible subsurface
$\Sigma_{0,4}$. These incompressible subsurfaces are unique up to 
homeomorphisms of the surface. Thus we may draw $(\alpha_1, \alpha_2)$
as in figure 4.
Then  $\phi(\alpha_i)$'s satisfy the  conditions that $\phi(\alpha_i)
\cap \phi(\alpha_{i+1}) \neq \emptyset$ and $\phi(\alpha_{i})
\cap \phi(\alpha_{i+2}) = \emptyset$. Now $\phi(\alpha_1) \perp_0
\phi(\alpha_2)$ follows from the lemma below.

{\bf Lemma 4.2.} \it Suppose $\beta_1, ..., \beta_5$ are five
pairwise distinct elements in $\Cal S(\Sigma_{0,5})$ so that
$\beta_i \cap \beta_{i+1} \neq \emptyset$ and $\beta_i \cap \beta_{i+2}
= \emptyset$ for all indices $i$ (mod 5). Then $\beta_i \perp_0 \beta_{i+1}$
for all $i$. \rm

{\it Proof.} We shall prove $\beta_1 \perp_0 \beta_2$
only. Take $b_i \in \beta_i$ so that $|b_i \cap b_j| = I(b_i, b_j)$.
Consider the subsurface $\Sigma_{0,4}$ bounded by $b_4$. The subsurface $\Sigma_{0,4}$ contains $b_1$ and $b_2$ by the assumption. Since $b_1 \cap b_3 =
\emptyset$, we conclude that $b_3 \cap \Sigma_{0,4}$ consists of parallel
copies of an arc in $\Sigma_{0,4}$. Furthermore, $b_1$ is determined
up to isotopy by $b_3$ and $b_4$. Indeed,  $b_1$ is isotopic to a boundary 
component of $N( b_3 \cup b_4)$.  Another way to see it is to use the
following lemma.

{\bf Lemma 4.3.} \it  Given two distinct classes in 
$\Cal S(\Sigma_{0,5})$ (resp. in $\Cal S(\Sigma_{1,2})$),  there 
is at most one class in  $\Cal S(\Sigma_{0,5})$ (resp. in
$\Cal S(\Sigma_{1,2})$) which is disjoint from both classes.
\rm

To see the proof, we note  that  the  only incompressible subsurfaces 
of negative Euler number in the surface are $\Sigma_{0,3}$, $\Sigma_{0,4}$
and $\Sigma_{1,1}$. Lemma 4.3 follows by considering the
smallest subsurface containing  the given classes.

Back to the proof of lemma 4.2, we have the same conclusion that
$b_5 \cap \Sigma_{0,4}$ consists of parallel copies of an arc in
$\Sigma_{0,4}$ and $b_5$
is determined uniquely up to isotopy by $b_2$ and $b_4$. Since $b_3
\cap b_5 = \emptyset$, $b_1 \cap b_2$ consists of two points. This
shows that $\beta_1 \perp_0 \beta_2$. 
$\square$.

{\it Case 2}. The surface is $\Sigma_{1,2}$. Take 
 $\alpha_1 \top \alpha_2$. We shall discuss three subcases: (2.1) 
$\alpha_1 \perp \alpha_2$, (2.2) $\alpha_1 \perp_0 \alpha_2$ so that
one of $\alpha_i$ is separating, (2.3) $\alpha_1 \perp_0 \alpha_2$
so that both elements $\alpha_i$ are non-separating.

{\it Subcase 2.1.} If $\alpha_1 \perp \alpha_2$, we extend 
$\{\alpha_1, \alpha_2$\} to a ``pentagon"
set $\{\alpha_1, ..., \alpha_5\}$ as in figure 5(a) where
$\alpha_i \cap \alpha_{i+2} = \emptyset$, $\alpha_1 \perp_0 \alpha_5,$
 $\alpha_2 \perp \alpha_3$, $\alpha_3 \perp_0 \alpha_4$, and
$\alpha_4 \cap \alpha_5  \neq  \emptyset$. Now $\phi(\alpha_1) \perp
\phi(\alpha_2)$ follows by the same argument as in case 1 (applied to
$\Sigma_{1,1}$ instead of $\Sigma_{0,4}$). See figure 5(b).

\midspace{0.1cm}
\centerline{\epsfbox{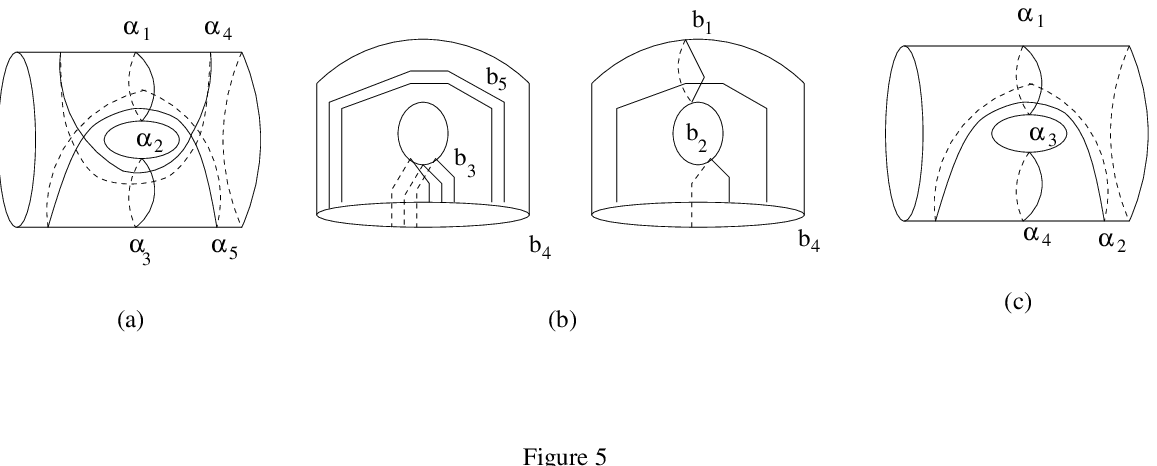}}
\midspace{0.1cm}

{\it Subcase 2.2}. If  $\alpha_1 \perp_0 \alpha_2$ so that $\alpha_2$
is separating, then $\alpha_1$ is non-separating. We extend 
it to a four-element set $\{\alpha_1, ..., \alpha_4\}$ as in figure 5(c) where
$\alpha_3 \cap \alpha_2 = \alpha_1 \cap \alpha_4 = \alpha_4 \cap \alpha_2 = \emptyset $ and $\alpha_1 \perp \alpha_3 \perp \alpha_4$. 
By subcase 2.1, we conclude that $\phi(\alpha_1) \perp \phi(\alpha_3)$
and $\phi(\alpha_3) \perp \phi(\alpha_4)$. 
Furthermore, $\phi(\alpha_2) \cap \phi(\alpha_3) = \phi(\alpha_2)
\cap \phi(\alpha_4) = \phi(\alpha_1) \cap \phi(\alpha_4) = \emptyset$. 
Now by lemma 4.3,  $\phi(\alpha_2)$  is
determined by $\phi(\alpha_3)$ and $\phi(\alpha_4)$.
Thus $\phi(\alpha_1) \perp_0 \phi(\alpha_2)$.

{\it Subcase 2.3.} 
If $\alpha_1 \perp_0 \alpha_2$  so that both  $\alpha_i$'s are 
non-separating, then both $\alpha_1 \alpha_2$ and
$\alpha_2 \alpha_1$ are separating. Since  $\alpha_1 \alpha_2
\perp_0 \alpha_i$ for $i=1,2$, by subcase 2.2, we obtain 
$\phi(\alpha_1 \alpha_2) \perp_0 \phi(\alpha_i)$ for $i=1,2$.
Similarly, $\phi(\alpha_2 \alpha_1) \perp_0 \phi(\alpha_i)$ for
$i=1,2$. Since $\alpha_1$, $\alpha_2$,  $\alpha_1 \alpha_2$
and $\alpha_2 \alpha_1$ are in a subsurface homeomorphic
to $\Sigma_{0,4}$,
by lemma 2.2, we conclude that classes 
$\phi(\alpha_1)$, $\phi(\alpha_2)$,  $\phi(\alpha_1 \alpha_2)$ 
and $\phi(\alpha_2 \alpha_1)$
are in a subsurface homeomorphic to $\Sigma_{0,4}$ as well.
Thus by
lemma 3.1(b) applied to the subsurface $\Sigma_{0,4}$, we have
$\phi(\alpha_1) \perp_0 \phi(\alpha_2)$.

To show  the last assertion in the proposition for $\Sigma_{1,2}$,
take $\alpha_1 \top \alpha_2$. Then 
$\alpha_1 \alpha_2$ is not $\top$-related to $\alpha_2 \alpha_1$ and
$\alpha_1, \alpha_2, \alpha_1 \alpha_2$ and $\alpha_2 \alpha_1$ 
are in a subsurface homeomorphic to $\Sigma_{1,1}$ or
$\Sigma_{0,4}$. Since $\phi$ preserves disjointness and relations
$\perp$ and $\perp_0$,  $\phi(\alpha_1)$, $\phi(\alpha_2)$, 
$\phi(\alpha_1 \alpha_2)$ and $\phi(\alpha_2 \alpha_1)$ are in a 
subsurface homeomorphic to $\Sigma_{1,1}$ or $\Sigma_{0,4}$ and
$\phi(\alpha_1 \alpha_2)$ is not $\top$-related to $\phi(\alpha_2
\alpha_1)$.
Applying lemma 3.1 (b) to the subsurface,
we conclude that 
$\{ \phi(\alpha_1  \alpha_2), \phi(\alpha_2 \alpha_1)\}
= \{ \phi(\alpha_1) \phi(\alpha_2), \phi(\alpha_2) \phi(\alpha_1)\}$.

We now prove the proposition by induction on  $|\Sigma_{g,n}| = 
3g +n$. The result holds for $|\Sigma| =5$ by the above two cases. 
If  $|\Sigma| \geq 6$,  take  $\alpha \top \beta$ in $\Cal S(\Sigma)$. 
Then $\alpha$ and  $\beta$  lie in an incompressible subsurface 
homeomorphic to $\Sigma_{1,1}$ or $\Sigma_{0,4}$.  
Choose a class $\gamma$ disjoint from $\alpha$ and $\beta$ so that 
either $\gamma$  is non-separating or is a boundary class. Take
$c \in\gamma$  and let $\Sigma'$ be a component of $\Sigma - int(N(c))$
 which contains $\alpha$ and $\beta$.  Then  $|\Sigma'| \geq 5$ by the
choice of $\gamma$.  Since $|\Sigma| \geq 6$, by lemma 2.3, there is 
a homeomorphism $h$ of the surface sending $\gamma$ to $\phi(\gamma)$.
After composing $\phi$ by $h^{-1}$, we may assume that 
 $\phi(\gamma) =\gamma$. It follows that $\phi(\Cal S(\Sigma')) =
\Cal S(\Sigma')$ by the choice of $\gamma$. 
Thus by lemma 2.3, $\phi|_{\Cal S(\Sigma')}$ preserves the separating 
classes if  $\Sigma' \ncong \Sigma_{1,2}$. If $\Sigma' \cong \Sigma_{1,2}$,
then by corollary 2.4, $\phi|_{\Cal S(\Sigma')}$ again preserves the 
separating classes. Thus by the induction hypothesis applied to $\Sigma'$,
we conclude that if $\alpha \perp \beta$, then $\phi(\alpha) \perp \phi(\beta)$
and if $\alpha \perp_0 \beta$ then $\phi(\alpha) \perp_0 \phi(\beta)$.
 Furthermore, in both cases, we have
 $\{ \phi(\alpha \beta), \phi(\beta \alpha)\}
=\{ \phi(\alpha) \phi(\beta), \phi(\beta) \phi(\alpha)\}$.
$\square$

\S5.  {\bf Proof of  the Main  Theorem }

Recall that surfaces in this section have negative Euler number. 
By proposition 4.1 and lemma 2.2,  it suffices to show the following
in order to finishing proof of the main theorem. 

{\bf Theorem.} \it Suppose $\phi: \Cal S(\Sigma) \to \Cal S(\Sigma)$ is 
a bijection
preserving disjointness, the separating classes, the relations $\perp$,
$\perp_0$, and $\{\phi(\alpha \beta), $$\phi(\beta \alpha)\}
=\{ \phi(\alpha) \phi(\beta),$  \newline
 $ \phi(\beta) \phi(\alpha)\}$. Then $\phi =
h$ for some $h \in $ Home($\Sigma$). \rm

{\bf Proof.} We use induction on $|\Sigma|$. For $|\Sigma| =4$, the result 
follows from lemma 3.1. If $|\Sigma| \geq 5$, we decompose
$\Sigma = X \cup Y$ where $X$ and $Y$ are compact incompressible subsurfaces
so that the following conditions hold: 
(i) $X \cap Y \cong \Sigma_{0,3}$,
(ii) if the genus $g=0$, then $X \cong \Sigma_{0,4}$ and
 $Y \cong \Sigma_{0,n-1}$, 
(iii) if the genus $g \geq 1$, then 
$X \cong \Sigma_{1,1}$ and
$Y \cong \Sigma_{g-1, n+2}$.
See figure 6. 

\midspace{0.1cm}
\centerline{\epsfbox{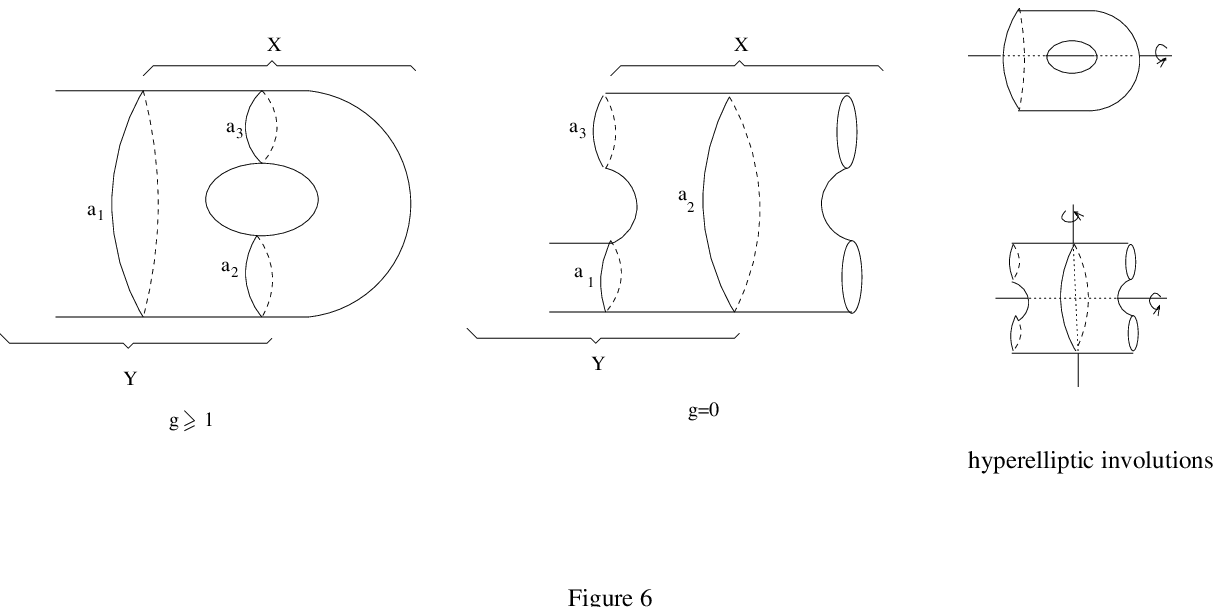}}
\midspace{0.1cm}

We write $\partial(X \cap Y) = a_1 \cup a_2 \cup a_3$ 
so that $a_1 \subset \partial X$, $a_2 \cup a_3 \subset \partial Y$,
and if the genus $g=0$, $a_3 \subset \partial \Sigma$.
By corollary 2.5, we find 
$h_1 \in$ Home($\Sigma$) so that $h_1( \phi([a_i])) = [a_i]$ for $i=1,2$. 
Thus, by replacing $\phi$ by $h_1 \phi$,  we may assume that
$\phi([a_i]) = [a_i]$ for $i=1,2$. 
This implies that $\phi(\Cal S(X)) = \Cal S(X)$
and $\phi(\Cal S(Y)) =  \Cal S(Y)$. 
Now by the construction, $|X|, |Y|$
$< |\Sigma|$ and $Y \ncong \Sigma_{0,3}$. 
We claim that the restrictions of $\phi$ to $\Cal S(X)$ and
$\Cal S(Y)$ satisfy the induction hypothesis. Evidently the restrictions
preserve the disjointness, the relations $\perp$ and $\perp_0$ and
 $\{\phi(\alpha \beta), \phi(\beta \alpha)\}
=\{ \phi(\alpha) \phi(\beta), \phi(\beta) \phi(\alpha)\}$. By lemma 2.3 
and corollary 2.4, the restriction of $\phi$ to $\Cal S(Y)$ preserves
the separating classes.
Thus, by the induction hypothesis, we find
$h_X \in $ Home($X$), $h_Y \in$ Home($Y$) so that $h_X = \phi|_{\Cal S(X)}$,
and $h_Y = \phi  |_{\Cal S(Y)}$.  

We shall use the following results to finish the proof of the
theorem. The proofs of these results are deferred to the end of
this section.

{\bf Lemma 5.2.} \it We may modify $h_X$  and $h_Y$ by composing with  
hyperelliptic involutions which are in the center of the mapping class group so
that after the modification $h_X(a_i) = h_Y(a_i)$, for $i=1,2,3.$  \rm

{\bf Proposition 5.3.} \it Both homeomorphisms $h_X$ and $h_Y$ are
orientation preserving or both are orientation reversing. \rm

{\bf Lemma 5.4.} \it An orientation preserving homeomorphism of
the 3-holed sphere leaving each boundary component invariant
is isotopic to the identity map. \rm

By lemmas 5.2, 5.4 and proposition 5.3, we conclude that
$h_X|_{X \cap Y} : X \cap Y  \to \Sigma$ and $h_Y|_{X \cap Y}: X \cap Y
 \to \Sigma$
are isotopic. Thus there exists $h
\in $ Home($\Sigma$) so that $h|_X \cong h_X$ and $h|_Y \cong h_Y$.
We have $\phi|_{\Cal S(X) \cup \Cal S(Y)} = h|_{\Cal S(X) \cup \Cal S(Y)}$.
The aim is to show that $\phi = h$. Since 
 $\{\phi(\alpha \beta), \phi(\beta \alpha)\}
=\{ \phi(\alpha) \phi(\beta), \phi(\beta) \phi(\alpha)\}$ and
 $\{h(\alpha \beta), h(\beta \alpha)\}
=\{ h(\alpha) h(\beta), h(\beta) h(\alpha)\}$, by corollary 3.3, it
suffices to show that $h(\alpha) = \phi(\alpha)$ for all
$\alpha$ so that $\alpha \perp_0 [a_1]$ and either $\alpha \top
[a_i]$ or $\alpha \cap [a_i] = \emptyset$ for $i=2,3$. 
Since $[a_1]$ is separating, $I(\alpha, a_2) + I(\alpha, a_3)$
is even. Thus $(I(\alpha, a_2), I(\alpha, a_3))$ is one of
the following four pairs $(0,2), (2,0), (1,1)$ or $(2,2)$. On the
other hand, $a_3$ is either a boundary component or is isotopic
to $a_2$ by the construction. Thus $(I(\alpha, a_2), I(\alpha, a_3))=
(0,2)$ is impossible.  We shall discuss the 
three cases $(I(\alpha, a_2), I(\alpha, a_3))= (0,2), (2,0)$ and $(1,1)$
separately.

 The strategy to show $h(\alpha) = \phi(\alpha)$ for these specific 
elements  $\alpha$  is as follows. First we construct an incompressible 
subsurface  $\Sigma' \cong \Sigma_{0,5}$ or $\Sigma_{1,2}$ which
contains both  $X$ and $\alpha$.  Second, we shall construct 
two distinct elements $[b_1]$ and $[b_2]$ in $(\Cal S(X) \cup \Cal S(Y))
\cap \Cal S(\Sigma')$  which are disjoint from $\alpha$. 
By the assumption,  $h([b_i]) = \phi([b_i])$ for $i=1,2$
and each elements of $\{ \phi(\alpha), h(\alpha)\}$ is disjoint from 
 $\{ h([b_1]), h([b_2])\}$.  Finally, we show that $\phi(\alpha)$
is in  $\Cal S(h(\Sigma'))$. By lemma 4.3 applied to $h(\Sigma')$ 
 and the pair $\{ h([b_1]), h([b_2])\}$,  we conclude that 
$h(\alpha) = \phi(\alpha)$.

Now take $s \in \alpha$ so that $|s \cap a_i | = I(\alpha , a_i)$.

{\it  Case 1.} $(I(\alpha, a_2), I(\alpha, a_3)) =(0,2)$ and
$\alpha \perp [a_2]$. Then the surface $X \cong \Sigma_{0,4}$.

\midspace{0.1cm}
\centerline{\epsfbox{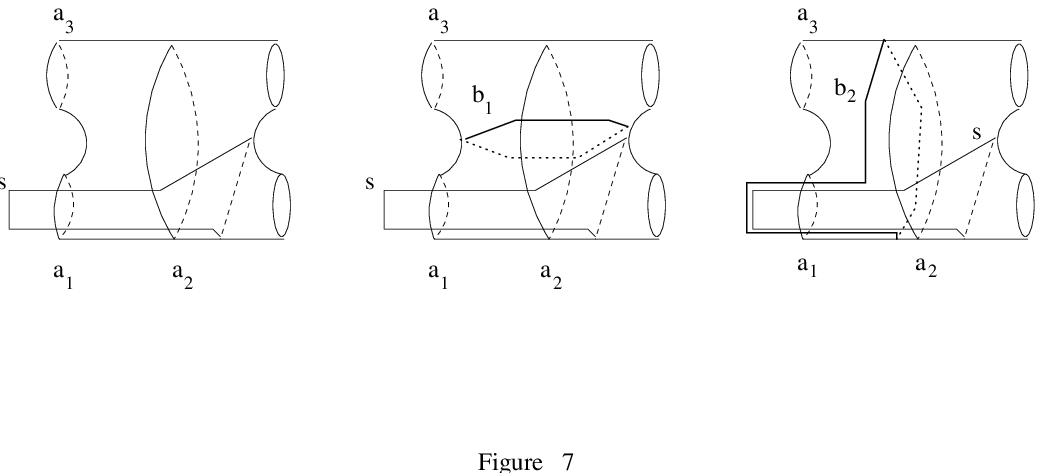}}
\midspace{0.1cm}

Let $\Sigma' = N(s) \cup X \cong \Sigma_{0,5}$. Then
$\Sigma'$ is incompressible. Take two essential
nonboundary parallel simple loops $b_1$ and $b_2$ in $\Sigma'$
so that (i) $b_1 \subset X$ and $b_2 \subset Y$, (ii)
$[b_i] \cap \alpha = \emptyset$, for  $i = 1, 2$
and (3) $ [b_1] \neq [b_2]$ as in figure 7.

The isotopy classes of each boundary component of $\Sigma'$ is 
either in $\partial \Sigma$ or is 
in $\Cal S(X) \cup \Cal S(Y)$. By the assumption, we have
$h(\beta)  = \phi(\beta)$ for each isotopy class $\beta$ of the
component of $\partial \Sigma'$ so that $\beta \in \Cal S(\Sigma)$.
Now  $\phi(\alpha)$ is disjoint from the isotopy classes 
of the boundary components of $h(\Sigma')$ and
$\phi(\alpha)$ intersets an isotopy class in $\Cal S(h(\Sigma'))$.
This shows that $\phi(\alpha)$ is in $\Cal S(h(\Sigma'))$.
Furthermore,  $\phi(\alpha)$ and $h(\alpha)$ 
are disjoint from  $h([b_i])$$ ( = \phi([b_i]))$ for $i=1,2$. Thus 
by lemma 4.3 applied to $h(\Sigma')$, $\phi(\alpha) = h(\alpha)$.

{\it Case 2}.
$(I(\alpha, a_2), I(\alpha, a_3)) =(1,1)$. Then $X \cong \Sigma_{1,1}$. Let
$\Sigma' = N(s) \cup X$. Then $\Sigma'$ is incompressible and is
homeomorphic to $\Sigma_{1,2}$. Choose two non-isotopic, non-boundary
parallel curves $b_1$ and $b_2$ in $\Sigma'$ as in figure 8. By the
construction, $[b_i] \in \Cal S(X) \cup \Cal S(Y)$ and $[b_i] \cap \alpha
= \emptyset$ for $i=1,2$. Furthermore, each component of
$\partial \Sigma'$ is either in $X, Y$ or in $\partial \Sigma$. Thus
$\phi(\alpha) = h(\alpha)$ by the same argument as in case 1.

\midspace{0.1cm}
\centerline{\epsfbox{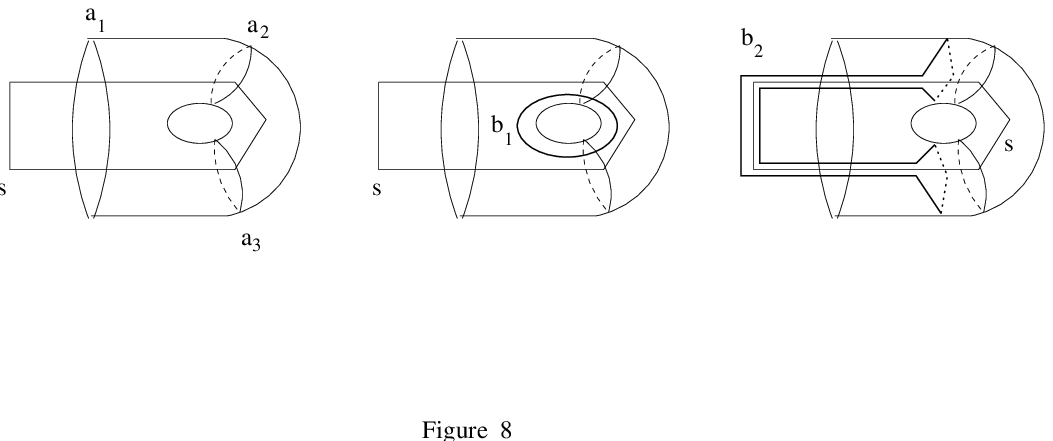}}
\midspace{0.1cm}

{\it Case 3}. $(I(\alpha, a_2), I(\alpha, a_3) ) =(2,2)$ and
$\alpha \perp_0 [a_2]$. This case does not occur. Indeed, $\alpha \perp_0
[a_2]$ shows that the arc $s \cap X$ would intersect $a_2$ at two
points of different signs. Since $X \cong \Sigma_{1,1}$, this shows
that $I(s, a_1) = 0$. This contradicts the assumption that $\alpha \perp_0
[a_1]$. $\square$

We now prove lemmas 5.2 and Proposition 5.3. Lemma 5.4 is well known. See 
for instance [FLP], expos\'e 2.

{\it Proof of lemma 5.2.} Since hyperelliptic involutions of $X$
act trivially on S(X), by composing $h_X$  by an isotopy and hyperelliptic 
involutions, we may assume that  $h_X(a_i) = a_i$ for $i=1,2, 3$.
Since $h_Y(a_1) \cong  a_1$,  we may assume that $h_Y(a_1) = a_1$ after
an isotopy.
        
If  $Y \cong \Sigma_{0,4}$,  then we may  assume that 
$ h_Y(a_i) = a_i$ for $i=2,3$ by composing $h_Y$ 
by hyperelliptic involutions. Thus the lemma follows in this case.

 If $Y \ncong \Sigma_{0,4}$, then $h_Y$ permutes $\{[a_2] , [a_3]\}$.
If $g \geq 1$,  by composing  $h_X$ by
hyperelliptic involutions if necessary, we obtain $h_X(a_i) = h_Y (a_i)$
for $i=1,2,3$. 
 If the genus $g=0$,  we shall prove that $h_Y(a_i) \cong a_i$ for
$i=2,3$. Suppose otherwise that  $h_Y(a_2) \cong a_3$.
Choose a boundary class  $\beta \in \Cal S(Y)$  so that 
$\beta$, $a_3$ and a component $b$ of $\partial Y \cap \partial \Sigma$
bound $\Sigma_{0,3}$ in $Y$. 
Thus  $\beta$ is also a boundary class
 in $\Sigma$. By lemma 2.3, $h_Y(\beta) $ $(= \phi(\beta))$ is again
a boundary class in $\Sigma$. But $h_Y(\beta)$ is also
a boundary class in $Y$ since $h_Y(\beta), a_2 =h_Y(a_3)$ and
$h_Y(b)$ bound a 3-holed sphere in Y. Since $[a_2] \in$ $\Cal S(\Sigma)$, this shows
that $Y \cong \Sigma_{0,4}$ which contradicts the assumption. $\square$

{\it 	Proof of proposition 5.3.}
Suppose otherwise, we may assume that $h_X$  
is orientation reversing and  $h_Y$ is orientation preserving. Thus
$\phi(\alpha) \phi(\beta) = \phi(\alpha \beta)$ for
$\alpha \top \beta$ in $ \Cal S(Y)$ and
$\phi(\beta) \phi(\alpha) = \phi(\alpha \beta)$
for $\alpha \top \beta$ for $\alpha, \beta$ in  $\Cal S(X)$.

If the genus $g=0$,  construct two curves $x$ and  $y$ as in  figure 9 
so that  $[x] \in \Cal S(X)$, $[y] \in \Cal S(Y)$, $[x] \perp_0 [y]$,
$[x] \perp_0 [a_2]$, $[y] \perp_0 [a_1]$, and $|y \cap a_1| =2$.
Then  $\phi(a_2 x) = \phi(x) \phi(a_2)$ and $\phi(a_1 y) =\phi(a_1) \phi(y)$.
Furthermore, the subsurface $\Sigma' = N(y) \cup X \cong
\Sigma_{0,5}$ is incompressible in  $\Sigma$. 
Thus two classes $\alpha, \beta \in \Cal S(\Sigma')$ are
 disjoint in $\Cal S(\Sigma)$ if and only if they are disjoint in  $\Cal S(\Sigma')$.
 We now use lemma 5.5 below to derive a contradiction. 

{\bf Lemma 5.5.} \it 
Suppose $\alpha \perp_0 \beta \perp_0 \gamma$, 
$\alpha \cap \gamma = \emptyset$ in $\Cal S(\Sigma_{0,5})$. Then 
$\alpha \beta \cap \gamma \beta = \emptyset$, $\alpha \beta \cap \beta \gamma
\neq \emptyset$, $\beta \alpha \cap \gamma \beta \neq \emptyset$, and $
\beta \alpha \cap \beta \gamma = \emptyset$. \rm

\midspace{0.1cm}
\centerline{\epsfbox{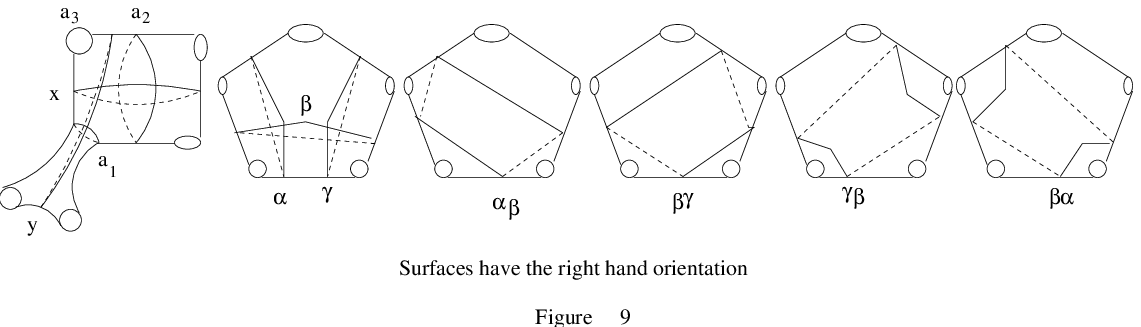}}
\midspace{0.1cm}

{\it Proof.}
Take a triple $(\alpha, \beta, \gamma)$ as in figure 9. Then
the lemma follows for the triple in figure 9 by the calculation
in the figure. On the other hand, there is only one triple 
$(\alpha, \beta, \gamma)$ 
satisfying the conditions in the   lemma up to self-homeomorphisms
of the surface.  Thus the lemma follows.
To see the uniqueness of the triple $(\alpha, \beta, \gamma)$ up to
homeomorphisms, we take three representatives $a,b, c$ in $\alpha,
\beta, \gamma$ respectively so that they intersect minimally. Then
the surface $\Sigma_{0,5}$ is homeomorphic to a regular neighborhood
$N(a \cup b \cup c)$. Furthermore, the union $a\cup b \cup c$  is unique
up to homeomorphisms. Thus the assertion follows. $\square$

Apply lemma 5.5 to $(\alpha, \beta, \gamma)$ = $([a_1], [y], [x])$
and $(\phi(a_1),$$ \phi(y),$$ \phi(x))$. We obtain

$$ [a_1][y] \cap [x][y] =\emptyset \quad [a_1][y] \cap [y][x] \neq \emptyset
\quad \tag 1$$

and 
$$ \phi(a_1) \phi(y) \cap \phi(x) \phi(y) = \emptyset \quad
\phi(a_1) \phi(y) \cap \phi(y) \phi(x) \neq \emptyset \tag 2$$

By applying $\phi$ to (1) and use $\phi(a_1 y) = \phi(a_1) \phi(y)$, we
obtain

$$ \phi(a_1) \phi(y) \cap \phi(xy) = \emptyset \quad \phi(a_1) \phi(y)
\cap \phi(yx) \neq \emptyset \tag 3$$

Since $\{\phi(xy), \phi(yx)\} =\{ \phi(x) \phi(y), \phi(y) \phi(x)\}$,
by comparing equations (2) and (3), we obtain

$$ \phi(xy) = \phi(x) \phi(y) \tag 4$$

If we apply lemma 5.5 to $(\alpha, \beta, \gamma)$ = $([a_2], $$
[x], [y])$ and $(\phi(a_2), $$\phi(x),$$ \phi(y))$ and
use $\phi(a_2x) = \phi(a) \phi(a_2)$, we obtain $\phi(xy) = \phi(y) \phi(x)$.
This contradicts (4).

If the genus $g=1$, we construct two curves $x, y$ as in figure 10
where $[x] \in \Cal S(X)$, $[y] \in \Cal S(Y)$ so that $[x] \perp [y]$,
$[x] \perp [a_2]$, $[y] \perp_0 [a_1]$ and $|y \cap a_1| =2$. The
subsurface $\Sigma' = N(y) \cup X \cong \Sigma_{1,2}$ is
incompressible in $\Sigma$. We use lemma 5.6 below to obtain a contradiction.

{\bf Lemma 5.6.} \it
If $\alpha \perp_0 \beta \perp \gamma \perp \delta$,
$\alpha \cap \delta = \beta \cap \delta = \beta \cap \gamma =\emptyset$ in
$\Cal S(\Sigma_{1,2})$, then $\alpha \beta \perp \gamma \beta$,
$\alpha \beta$ is not $\top$-related to  $\beta \gamma$, 
$\beta \alpha$ is not $\top$-related to $\gamma \beta$, and
$\beta \alpha \perp \beta \gamma$. Furthermore, $\delta \gamma \cap
\beta \gamma = \gamma \delta \cap \gamma \beta = \emptyset$, and
$\delta \gamma \cap \gamma \beta \neq \emptyset$, $\gamma \delta 
\cap \beta \gamma \neq \emptyset$. \rm

\midspace{0.1cm}
\centerline{\epsfbox{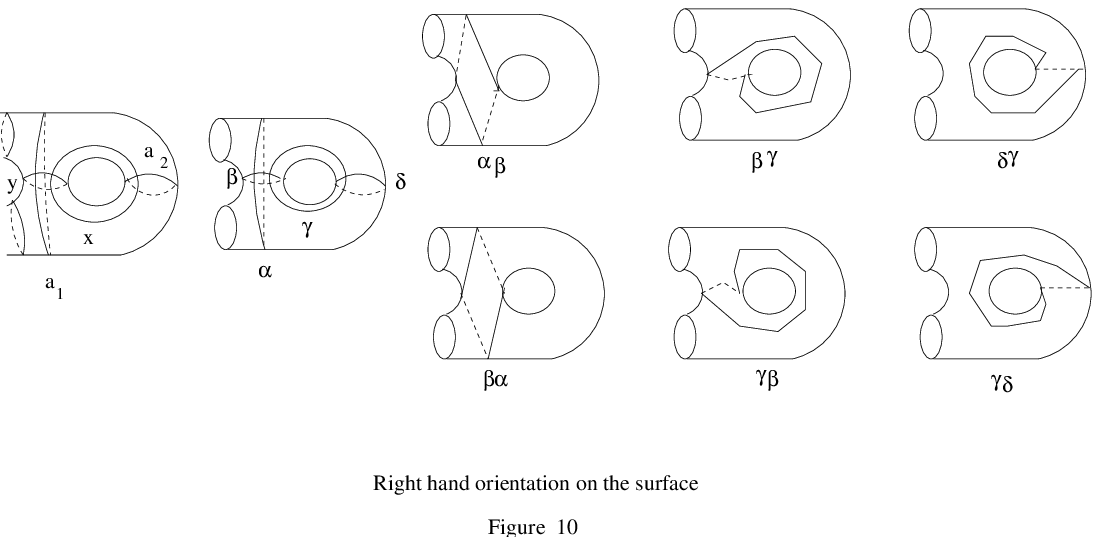}}
\midspace{0.1cm}

See figure 10 for a verification of the lemma for a
specific choice of the quadruple  $(\alpha, \beta, \gamma, \delta)$.
But the quadruple satisfying the conditions in the lemma is unique
up to self-homeomorphism of the surface. Indeed, by lemma 4.3,
$\delta$ is uniquely determined by $\alpha, \beta,$ and 
$\alpha$ is uniquely determined by $\delta, \gamma$. The
uniqueness of the triple ($\beta, \gamma, \delta)$
(resp. $(\alpha, \beta, \gamma)$)
follows by the same argument as in lemma 5.5.

Now the proof is similar to the previous case. Namely, by the choice
of $x,y$, we have $\phi(a_1y) =\phi(a_1) \phi(y)$ and $\phi(a_2 x)
= \phi(x) \phi(a_2)$. Applying lemma 5.6 to the triples
$([a_1], [x],[y])$ and $(\phi(a_1), \phi(y), \phi(x))$
(as $(\alpha, \beta, \gamma)$), we
obtain $\phi(xy) = \phi(x) \phi(y)$. If we apply the lemma to the
different triples $([a_2], [x], [y])$ and $(\phi(a_2), $\newline 
 $\phi(x),$ $ \phi(y))$ (as ($\beta, \gamma, \delta)$),
we obtain $\phi(xy) = \phi(y) \phi(x)$. This is a contradiction.
$\square$

{\it Acknowledgments.} I would like to thank X.S. Lin for many 
discussions, F. Bonahon for informing me Ivanov's work and N. Ivanov 
for sending his unpublished 1989 preprint. I thank
the referee for many helpful suggestions.
The work is supported in part by the NSF.

\centerline{\bf Reference}

[Bi] Birman, J.S., Braids, links, and mapping class groups. Ann. of Math. Stud., 82, Princeton Univ. Press,
 Princeton, NJ, 1975

[De] Dehn, M.: Papers on group theory and topology. J. Stillwell (eds.).
 Springer Verlag, Berlin-New York, 1987

[FLP] Fathi, A., Laudenbach, F., Poenaru, V.: Travaux de Thurston sur
 les surfaces. Asterisque 66-67, Societe Mathematique de France, 1979

[Gr] Grothendieck, A.: Esquisse d'un programme. in Geometric Galois action,
 London Math. Soc. Lecture Note Ser., 242, 5-48, Cambridge Univ. Press,
 Cambridge, 1997.

[Har] Harer, J.: Stability of the homology of the mapping class groups of
 orientable surfaces. Ann. of Math., 121 (2), (1985) 215-249

[Har1] Harer, J.: The virtual cohomological dimension of the mapping class
 group of an orientable surface. Invent. math. 84 (1986), 157-176

[Hay] Harvey, W.: Boundary structure of the modular group. In: Riemann
 surfaces and related topics: Proceedings of the 1978 Stony Brook
 Conference (State Univ. New York, Stony Brook, N.Y., 1978), 245-251,
 Ann. of Math. Stud., 97, Princeton Univ. Press, Princeton, N.J., 1981

[He] Hempel, J.: 3-manifolds from curve complex point of view, preprint

[Iv] Ivanov, N.V.: Automorphisms of complexes of curves and of Teichm\"uller
 spaces. Inter. Math. Res. Notice, 14, (1997), 651-666

[Iv1] Ivanov, N. V.: Complexes of curves and Teichm\"uller modular groups.
Uspekhi Mat. Nauk., 42 (3) (1987), 49-91

[Iv2] Ivanov, N.V.: Complexes of curves and Teichm\"uller spaces, Math.
 Notes, 49 (1991), 479-484

[Ko] Korkmaz, M: Complexes of curves and the mapping class groups, Top.
 \& Appl., to appear.

[Li] Lickorish, R.: A representation of oriented combinatorial 3-manifolds.
 Ann. Math., 72 (1962), 531-540

[Lu] Luo, F.: Geodesic length functions and Teichm\"uller spaces, 
J. Differential Geometry, 48 (1998), 275-317

[Lu1] Luo, F.: On non-separating simple closed curves in a compact surface.
 Topology, 36 (2), (1997), 381-410

[Lu2] Luo, F.: Grothendieck's reconstruction principle and 2-dimensional
topology and geometry, Communications in Contemporary Mathematics
to appear.

[MM] Masur, H., Minsky, Y.: Geometry of the complex of curves I:
 hyperbolicity, preprint 1996

[Pa] Patterson, D.B.: The Teichm\"uller spaces are distinct. Proc.
 Amer. Math. Soc., 35 (1972), 179-182; 38 (1973), 668-669

[Th] Thurston, W.: On the geometry and dynamics of diffeomorphisms of
 surfaces. Bull Amer. Math. Soc., 19 (2) (1988), 417-438

[Ti] Tits, J.: On buildings and their applications. In: James, R.D. (ed.)
 International Congress of Mathematicians, Vancouver, B.C., 1974, 221-228

[Vi] Viro, O.: Links, two-fold branched coverings and braids,
 Soviet Math. Sbornik, 87 no 2 (1972), 216-228

\end